\documentstyle[12pt, fullpage]{article}
\newcounter{ppp}
\newcounter{pdsix}

\newcommand{\iv}{^{-1}}

\newcommand{\pp}{{\cal P} }

\newcommand{\sss}{{\cal S} }

\begin{document}
     \title{Embeddings of relatively free groups into finitely presented
groups}
   \author{A.Yu.Ol'shanskii and M.V.Sapir}
     \maketitle

\newtheorem{theorem}{Theorem}[section]
\newtheorem{prop}[theorem]{Proposition}
\newtheorem{lm}[theorem]{Lemma}
\newtheorem{cor}[theorem]{Corollary}
\newtheorem{df}[theorem]{Definition}

     \section{Introduction}

     By the Higman embedding theorem \cite{Ro} every finitely generated
recursively presented group can be isomorphically embedded into
some finitely presented group. Moreover the embedding can preserve
the solvability of the word problem \cite{Claph} or even the
solvability of the conjugacy problem \cite{OScollins}. The
embedding can be quasi-isometric \cite{Ol1} and can at the same
time dramatically improve the isoperimetric function of the group
\cite{BORS}. Moreover if the word problem for a finitely generated
group $G$ is solvable in non-deterministic time $T(n)$ by a Turing
machine then $G$ is a quasi-isometric subgroup of a finitely
presented group $H$ whose Dehn function is polynomially equivalent
to $T(n)$ \cite{BORS}. Thus  the word problem of a finitely
generated group
 $G$ is in NP if and only if
$G$ is a quasi-isometric subgroup of a finitely presented group with polynomial
Dehn function  \cite{BORS}.

There exists (see, for example, Valiev \cite{Valiev}) a relatively
simple finite presentation of a group containing all recursively
presented groups, but the embedding of a concrete finitely
generated group in this universal group is not really effective.
The embeddings mentioned in the previous paragraph are not really
effective either. In both case one needs to know a precise
description of a Turing machine \cite{Ro} (or a Diophantine
equation \cite{Valiev, LS} or an $S$-machine \cite{BORS})
``solving" the word problem in the original group in order to find
the presentation of the bigger group or to find a copy of the
given group inside the universal group.

The main goal of our paper is to show that very often an embedding
of $G$ into $H$ can be given in a more explicit and
straightforward way. In particular, the description of the set of
defining words for $H$ is explicit. (More precisely the
$S$-machines used in these constructions are so small that it is
easy to write down all the commands, or simply forget about
machines altogether and write down all the relations of the
presentation of $H$.) We can also say a lot about the structure of
$H$, and about the way $G$ is embedded into $H$.

     In this paper, we consider two classes of groups. The first class
consists of relatively free groups of finite ranks in varieties of
groups. For example, we show that finitely generated free
solvable, or free Burnside groups of sufficiently large exponents
can be easily embedded into finitely presented groups with
polynomial isoperimetric functions. Notice that the isoperimetric
functions of these relatively free groups are not even defined
because these groups are (as a rule) infinitely presented. In
order to deal  with relatively free groups, we introduce the so
called {\em verbal isoperimetric function} of a relatively free
group. This function seems to be interesting in itself.

The second class of groups consists of one relator metabelian
Baumslag-Solitar groups $BS_{k,1}=\langle a,b\ |\ a^b=a^k\rangle$.
Here $x^y$ stands for $y\iv xy$. It is well known that the Dehn
function of the group $BS_{k,1}$ is exponential \cite{WPG}. It is
also well known that these groups are very ``stubborn": they
``resist" being embedded into groups with small isoperimetric
functions. For example, there is a conjecture that the Dehn
function of every 1-related group containing $BS_{k,1}$, $k\ge 2$
is exponential. In this paper, we present simple finite
presentations of groups with polynomial isoperimetric functions
which contain $BS_{k,1}$.

     Let us give the necessary definitions.
Let $H$ be a group given by a finite presentation $P=\langle
x_1,\dots,x_k\  |\ r_1,\dots, r_\ell\rangle$. A non-decreasing
function $f: {\bf N} \rightarrow {\bf N}$ is called an {\em
isoperimetric function} of this presentation if any word
$w=w(x_1,\dots,x_k)$ of length $\le n$, which is equal to $1$ in
$H$, can be written in the free group $F(x_1,\dots,x_k)$ as a
product of at most $f(n)$ conjugates of the relators of $H$. In
other words, $f(n)$ is an isoperimetric function of the
presentation $P$ if every loop of length $\le n$ in the Cayley
complex corresponding to the presentation $P$ has area at most
$f(n)$. Isoperimetric functions of different finite presentations
of the same group are {\em equivalent} in some natural sense
\cite{Gersten}, so one can talk about isoperimetric functions of
the group $H$ forgetting about its presentations. We do not
distinguish equivalent isoperimetric functions in this paper.
Following Gersten, the smallest isoperimetric function of a group
$H$ is called the {\em Dehn function} of $H$.

Now let us define the verbal isoperimetric  functions. It is easy
to see that any variety of groups ${\cal V}$ that can be defined
by a finite set of identities, can be also defined by a single law
$v=1$ for some word $v = v(x_1,\dots x_k)$ from the (absolutely)
free group $F$ of infinite rank. The verbal subgroup $V\le F$
consists of all words vanishing in all groups of the variety
${\cal V}$. Any word $w\in V$ is freely equal to a product
$$\prod_{i=1}^N u_i v(X_{i1},\dots,X_{im})^{\pm 1} u_i^{-1}\eqno
(1.1)$$ for some words $u_i$ and $X_{ij}$. We call a
non-decreasing function $f_v: {\bf N}\rightarrow {\bf N}$ a {\it
verbal isoperimetric function} of the word $v$, if for any word
$w\in V$ there is a representation (1.1), where
$\sum_{ij}|X_{ij}|\le f_v(|w|)$. The smallest verbal isoperimetric
function will be called the {\it verbal Dehn function}. Note that
the verbal Dehn function exists because in the definition of
verbal isoperimetric functions one may restrict oneself to words
in variables $x_1,..., x_n$ and there are only finitely many such
words of any given length.

     If an identity $v'(x_1,\dots,x_{m'})=1$ is equivalent to the
identity $v=1$ then every value of $v'$ is a product of a fixed
number of values of the word $v$, and vice versa. Thus the verbal
Dehn function of the variety $V$ does not depend on the choice of
a defining law if one identifies functions which are {\em
$\Theta$-equivalent}. Recall that two functions $f,g:{\bf
N}\rightarrow {\bf N}$ are $\Theta$-equivalent if there exist
positive constants $c,d$ such that for every $n\in {\bf N}$ we
have:

$$f(n)\le c g(n) \hbox{ and } g(n)\le d f(n).$$

     Assume that a word $w$ is a product of two words $w', w''$ which
depend on disjoint sets of variables, $S'$ and $S''$. Then a representation
$$w=\prod u_i v(X_{i1},\dots, X_{im})u_i^{-1}$$
induces similar representations of $w'$ and $w''$ if we substitute 1
for the variables from $S''$ or from $S'$. For the corresponding words
$X'_{ij}$ and $X''_{ij}$ one has the obvious inequality
$|X'_{ij}|+|X''_{ij}|\le |X_{ij}|$. Hence the following statement holds.

\begin{prop} \label{prop1}  Any verbal Dehn function $f$ is superadditive,
i.e. $f(n_1+n_2)\ge f(n_1)+f(n_2)$.
\end{prop}

\medskip

{\bf Remark.} Recall that the superadditivity of any (usual) Dehn
function is still an open problem (see Guba and Sapir
\cite{freeproducts}).
\medskip

     In this paper we produce explicite finite presentations of
groups $H(v,m)$ to prove the following

     \begin{theorem}\label{th1} Let $f(n)$ be the verbal isoperimetric
function of a group variety ${\cal V}$ defined by an identity
$v=1$. Then the free group $F_m({\cal V})$ of rank $m$ in the
variety ${\cal V}$ can be isomorphically embedded into a finitely
presented group $H = H(v,m)$ with an isoperimetric function
$n^2 f(n^2)^2$. The presentation of $H$ is explicitely constructed
given the word $v$, and the embedding is quasi-isometric.
\end{theorem}

Recall that a finitely generated subgroup $G$ of a finitely
generated group $H$ is {\em quasi-isometric} (in other
terminology, undistorted) if there is a positive constant $C$ such
that $|g|_G\le C |g|_H$ for any $g\in G$. Here $|g|_G$, $|g|_H$
are the lengths of $g$ in the fixed finite sets of generators of
$G$ and $H$ respectively. This notion is well defined, but the
constant $C$ does depend on the choice of the finite generating
sets.

\medskip
{\bf Remark.} Without the word ``explicit" and with slightly worse
estimate for the isoperimetric function of the bigger group this
theorem is a corollary of the main result of \cite{BORS}. Indeed,
it is easy to see that if $f$ is a verbal isoperimetric function
of a group variety ${\cal V}$ then the word problem in any
relatively free group in this variety can be solved by a
non-deterministic Turing machine in time $f(n)$. Thus by the main
result in \cite{BORS}, this group can be embedded into a finitely
presented group with isoperimetric function $n^2f(n^2)^4$ (compare
with the estimate $n^2f(n^2)^2$ in our theorem).
\medskip

For example, let ${\cal V}={\cal B}_n$ be the Burnside variety
consisting of all groups satisfying the identity $x^n=1$ for a
fixed large odd number $n$. Then the function $f(s)=s^4$ is a
verbal isoperimetric function for ${\cal V}$. One can refer, for
instance, to Storozhev's argument in \cite{Ol89}, $\S 28.2$: a
word $w$ of length $s$ that is equal to 1 in the free Burnside
group $B(m,n)$, is a product of $n$-th powers of some words $A_i$
where the lengths of $A_i$ are bounded by a linear function of
$s$, and the number of factors is bounded by a cubic function.
Thus the function $s^{18}$ is an isoperimetric function of
$H(v,m)$ in this case. It can be proved by refining Storozhev's
argument that the function $s^{1+\epsilon(n)}$ is a verbal
isoperimetric function for ${\cal B}_n$ where $\epsilon
(n)\rightarrow 0$ for $n\rightarrow\infty$ (R. Mikhajlov,
unpublished). Therefore for very large odd $n$ the group
$H(x^n,m)$ satisfies a verbal isoperimetric inequality with
$f(s)=s^{8+\epsilon}$ for a small $\epsilon >0$.

The Schreier -- Reidemeister rewriting shows by induction, that
the variety of all solvable groups of derived length at most $d$
has a polynomial verbal isoperimetric function. It is interesting
to calculate (up to the equivalence) the verbal Dehn function of
these varieties, or the varieties of all nilpotent groups of given
nilpotency classes, etc. Results of Yu. G. Kleiman \cite{Kleiman}
show that there exists an identity defining a solvable group
variety with non-recursive verbal Dehn function. Indeed, Kleiman
\cite{Kleiman} constructed a finitely based solvable variety with
undecidable identity problem: given a word $w$, it is impossible
to decide whether the law $w=1$ follows from the defining laws of
his variety. Clearly this implies that Kleiman's variety cannot
have a recursive verbal isoperimetric function.

Similar questions can be raised for the function which counts just
the minimal number of factors $N$ in (1.1). Notice that from the
decidability of the Diophantine theory of the free group (Makanin
\cite{Makanin}), it follows that this function is recursive if and
only if the corresponding verbal Dehn function is recursive.

\begin{theorem}\label{th2} The group $BS_{k,1}=\langle a,b\ |\ b^{-1}ab=a^k\rangle$
is quasi-isometrically embeddable into a finitely presented group $H_{k,1}$
with isoperimetric function $n^{10}$.
\end{theorem}

This theorem gives an upper bound for the Dehn function of an
appropriate finitely presented group containing $BS_{k,1}$ as a
subgroup, and in contrast to \cite{BORS},
the presentation of $H_{k,1}$
is explicit and $S$-machines are used only in an implicite way.

\section{Defining relations}

Let ${\cal V}$ be the variety of groups defined by an identity
$v(x_1,...,x_k)=1$ where $v$ is a cyclically reduced word, that
can be written as $y_1\dots y_{M-1}$ with $y_s \equiv x_{i_s}^{\pm 1}$.

First we define an auxiliary group $G(v,m)$ which plays the same
role as the group $G_N(\cal S)$ in \cite{BORS} (and bares some
similarity with the Boone group from \cite{Ro}). The set of
generators of $G(v,m)$ consists of the letters $a_1,\dots,a_m$;
$q_1,\dots,q_{M}$;  $k_1,\dots,k_N $ (where $N\ge 29$);
$r_{1,1},...,r_{m,k}$.

We present the set of relations of the group $G(v,m)$ in three
different ways. First we just list all the relations. Then we will
present these relations in terms of an $S$-machine, and finally we
shall present a picture which will contain all the relations, and
also show the main corollaries of the relations.

{\bf 2.A. The list of relations.}
The finite set of defining relators for $G(v,m)$ is given below by equalities
(2.1) - (2.6).

$$r_i^{-1} q_{s+1} r_i = a_j q_{s+1}  \eqno (2.1)$$
if  $i = (j,l)$ and $y_s\equiv x_\ell$;

$$ r_i^{-1} q_s r_i = q_s a_j^{-1} \eqno (2.2)$$
if  $i = (j,l)$ and $y_s\equiv x_\ell^{-1}$;

$$ r_i^{-1} q_s r_i = q_s \eqno (2.3)$$
for all combinations of $i$ and $s$ which are not considered in (2.1) or (2.2);

$$ r_i^{-1} a_j r_i = a_j,\; i\in I, 1\le j\le m ;\eqno (2.4)$$


$$ r_i^{-1} k_j r_j = k_j,\; i\in I, 1\le j\le N; \eqno (2.5)$$

$$ k_1(q_1 q_2\dots q_M)\dots k_N(q_1 q_2 \dots q_M)=1.\eqno (2.6)$$

\medskip

{\bf 2.B. $S$-machines.} Let us give a precise definition of
$S$-machines \cite{SBR}. Let $k$ be a natural number.  Consider
now a {\em language of admissible words}. It consists of words of
the form $$q_1u_1q_2...u_kq_{k+1}$$ where $q_i$ are letters from
disjoint sets $Q_i$, $i=1,...,k+1$, $u_i$ are reduced group words
in an alphabet $Y_i$ ($Y_i$ are not necessarily disjoint), the
sets $\bar Y=\bigcup Y_i$ and $\bar Q=\bigcup Q_i$ are disjoint.

Notice that in every admissible word, there is exactly one
representative of each $Q_i$ and these representatives appear in
this word in the order of the indices of $Q_i$.

If $0\le i\le j\le k$ and $W=
q_1u_1q_2...u_kq_{k+1}$ is an admissible word then
the subword $q_iu_i...q_j$ of $W$ is called the $(Q_i,Q_j)$-subword of $W$
($i<j$).

An $S$-machine is a rewriting system \cite{KharSap}. The objects
of this rewriting system are all admissible words.

The rewriting rules, or {\em $S$-rules}, have the following form:
$$[U_1\to V_1,...,U_m\to V_m]$$ where the following conditions hold:

\begin{description}
\item Each $U_i$ is a subword of an admissible word starting with
a $Q_\ell$-letter and ending with a $Q_r$-letter (where $\ell=\ell(i)$
must not exceed $r=r(i)$, of
course).
\item If $i<j$ then $r(i) < \ell(j)$.
\item Each $V_i$ is also a subword of an admissible word whose $Q$-letters
belong to $Q_{\ell(i)}\cup...\cup Q_{r(i)}$ and which contains a $Q_\ell$-letter
and a $Q_r$-letter.
\item If $\ell(1)=1$ then $V_1$ must start with a $Q_1$-letter and if
$r(m)=k+1$ then $V_n$ must end with a $Q_{k+1}$-letter (so tape letters
are not inserted to the left of $Q_1$-letters and to the right of $Q_{k+1}$-letters).
\end{description}

To apply an $S$-rule to a word $W$ means to replace simultaneously subwords
$U_i$ by subwords $V_i$, $i=1,...,m$. In particular, this means that our
rule is not applicable if one of the $U_i$'s is not a subword of $W$. The
following convention is important:

After every application of a rewriting rule, the word is automatically
reduced. We do not consider reducing of an admissible word a separate step
of an $S$-machine.

We also always assume that an $S$-machine is symmetric, that is for every
rule of the $S$-machine the inverse rule (defined in the natural way) is
also a rule of this $S$-machine.

Notice that virtually any $S$-machine is highly nondeterministic.

Among all admissible words of an $S$-machine we fix one word
$W_0$. If an $S$-machine $\sss$ can take an admissible word $W$ to
$W_0$ then we say that $\sss$ {\em accepts} $W$. We can define the
time function of an $S$-machine as usual. If $Z\to Z_1\to...\to
Z_n=W_0$ is an accepting computation of the $\sss$-machine $\sss$
then $|Z|+|Z_1|+...+|Z_n|$ is called the {\em area} of this
computation. This allows us to define the {\em area function} of
an $S$-machine.

In \cite{SBR}, it is showed how to associate a group with any $S$-machine.
Our group $G(v,m)$ is the group associated with the following simple $S$-machine
${\cal S}(v)$.

Its language of admissible words coincides with the set of words
of the form $q_1u_1q_2u_2...q_M$ where $u_i$ are any group words
in the alphabet $\{a_1,...,a_m\}$. Each command of this
$S$-machine corresponds to a variable of $v$ and a letter from
$\{a_1,...,a_m\}$. Let $x$ be one of the variables of $v$, which
occurs with exponent $+1$ at positions $i_1,...,i_s$ of $v$ and
occurs with exponent $-1$ at positions $j_1,...,j_t$ and let $a$
be any letter from $\{a_1,...,a_m\}$. Then the corresponding
command multiplies $q_{i_1+1},..., q_{i_s+1}$ by $a$ on the left,
multiplies $q_{j_1},...,q_{j_t}$ by $a^{-1}$ on the right, and
does not change other $q$.

For example, if $v=x^{M-1}$ then $v$ contains only one variable,
and commands of the $S$-machine ${\cal S}(v)$ are indexed by
letters from $\{a_1,...,a_m\}$ and each command has the form
$[q_2\to a_iq_2,...,q_M\to a_iq_M]$.

Let $W_0=q_1...q_M$. Here is the main property of the $S$-machine
${\cal S}(v)$. This statement can be easily proved by induction.

{\bf The main property of $\sss(v)$.} An admisible word
$q_1u_1q_2...u_{M-1}q_M$ is accepted by the $S$-machine ${\cal
S}(v)$ if and only if $u_i\equiv u_j \iff y_i=y_j$. In other words
acceptable words are obtained from the words of the form
$v(u_1,...,u_k)=u_{i_1}\cdot...\cdot u_{i_{M-1}}$ by inserting
$q_1,...,q_M$ between the factors $u_{i_j}$.

For example, if $v=x^{M-1}$ then accepted words have the form
$q_1uq_2u...uq_{M}$ where $u$ is an arbitrary word in the alphabet
$\{a_1,...,a_n\}$.
\medskip

The construction (essentialy from \cite{SBR}) which simulates an
$S$-machine in a group is the following.

Let $\Theta$ be the set of rules of $\sss$. Let us call one of each pair
of mutually inverse rules from $\Theta$ {\em positive} and the other
one {\em negative}. The set of all positive rules will be denoted
by $\Theta_+$.

Let $N$ be any positive integer ($\ge 29$). Let $A$ be the set of all letters occurring in the admissible words of $\sss$ union with the set $\{k_j\ |\ j=1,\ldots, N\}\cup
\Theta_+.$

Our group is generated by the set $A$ subject to the set
$\pp_N(\sss)$ of relations described below.

{\bf 1. Transition relations}. These relations correspond to
elements of $\Theta_+$.

Let $r\in \Theta_+$, $r=[U_1\to V_1,...,U_p\to V_p]$. Then we
include relations $U_1^r=V_1,...,U_p^r=V_p$ into $\pp_N(\sss)$. If
some set $Q_j$ does not have a representative in any of the words
$U_i$ then we include all the commutativity relations $q^r=q$,
$q\in Q_j$.

{\bf 2. Auxiliary relations}.

These are all possible relations of the form $r x=xr$ where
$x$ is one of the letters in $\{a_1,...,a_m,k_1,...,k_N\}$, $r\in \Theta^+$.

{\bf 3. The hub relation.}

$$k_1W_0k_2W_0...k_NW_0=1$$

It is easy to see that the group presentation associated with our
$S$-machine $S(v)$ coincides with the presentation constructed
above in section A.

{\bf 2.C. A picture.}  For simplicity, let us take $v=x^n$: our
construction does not depend much on the word $v$ anyway.

Further simplifying the situation (and the future picture) let us
take $n=3$. The construction really does not depend much on $n$
either, so we shall sometimes write $n$ instead of $3$.
\setcounter{pdsix}{\value{ppp}} Figure \thepdsix\  shows the van
Kampen diagram (below it will be called a {\em disc}) with
boundary label
$\Sigma(q_1uq_2uq_3uq_4)=k_1q_1uq_2uq_3uq_4k_2q_1uq_2uq_3uq_4k_3...k_{N}q_1uq_2uq_3uq_4$.

\bigskip

\unitlength=0.8mm
\linethickness{0.4pt}
\begin{picture}(140.00,92.00)
\put(89.17,46.00){\oval(101.00,83.33)[]}
\put(88.33,45.50){\oval(86.00,69.67)[]}
\put(88.67,44.50){\oval(73.33,56.33)[]}
\put(88.83,44.17){\oval(59.67,44.33)[]}
\put(89.33,43.50){\oval(48.67,32.33)[]}
\put(89.33,43.00){\oval(37.33,22.67)[]}
\put(90.67,42.83){\oval(25.33,13.00)[]}
\put(78.00,43.00){\line(-1,1){39.33}}
\put(44.00,87.67){\line(1,-1){38.33}}
\put(97.00,49.33){\line(1,1){38.00}}
\put(103.00,44.67){\line(1,1){36.67}}
\put(38.67,10.00){\line(4,3){40.33}}
\put(44.33,4.33){\line(6,5){39.00}}
\put(133.67,4.33){\line(-1,1){33.67}}
\put(100.00,38.00){\line(0,0){0.00}}
\put(139.67,10.00){\line(-6,5){37.00}}
\put(60.00,72.67){\vector(1,0){3.33}}
\put(69.00,72.67){\vector(1,0){5.33}}
\put(97.33,72.67){\vector(1,0){6.00}}
\put(114.33,72.67){\vector(1,0){6.00}}
\put(68.00,72.67){\circle*{0.67}}
\put(96.33,72.67){\circle*{0.67}}
\put(113.67,72.67){\circle*{0.67}}
\put(55.00,80.33){\vector(1,0){3.33}}
\put(69.00,80.33){\vector(1,0){5.33}}
\put(97.33,80.33){\vector(1,0){6.00}}
\put(68.00,80.33){\circle*{0.67}}
\put(96.33,80.33){\circle*{0.67}}
\put(120.67,80.33){\circle*{0.67}}
\put(68.00,72.67){\line(-5,6){6.33}}
\put(74.33,72.67){\line(0,1){7.67}}
\put(96.33,72.67){\line(-4,5){6.00}}
\put(64.33,72.67){\line(-4,5){6.33}} 
\put(103.33,72.67){\line(0,1){7.67}}
\put(39.33,87.33){\makebox(0,0)[cc]{$k_1$}}
\put(140.00,86.67){\makebox(0,0)[cc]{$k_2$}}
\put(139.67,3.33){\makebox(0,0)[cc]{$k_3$}}
\put(65.00,70.33){\makebox(0,0)[cc]{$q_1$}}
\put(71.00,70.33){\makebox(0,0)[cc]{$q_2$}}
\put(99.00,70.67){\makebox(0,0)[cc]{$q_3$}}
\put(55.00,82.33){\makebox(0,0)[cc]{$q_1$}}
\put(67.33,82.00){\makebox(0,0)[cc]{$aq_2$}}
\put(97.00,82.67){\makebox(0,0)[cc]{$aq_3$}}
\put(117.00,70.67){\makebox(0,0)[cc]{$q_4$}}
\put(120.67,82.67){\makebox(0,0)[cc]{$aq_4$}}
\put(82.33,92.00){\makebox(0,0)[cc]{$q_1uq_2uq_3uq_4$}}
\put(113.67,72.67){\line(0,1){7.67}}
\put(120.67,80.33){\vector(1,0){3.67}}
\end{picture}
\begin{center}
\nopagebreak[4]
Fig. \theppp.

\end{center}
\addtocounter{ppp}{1}

On the boundary of this diagram we can read the word
$\Sigma(q_1uq_2uq_3uq_4)$. The
words on each of the concentric circles is labeled by $\Sigma(q_1u_iq_2u_iq_3u_iq_4)$
where $u_i$ is a prefix of $u$ of length $i-1$. The word written on the
innermost circle is the hub, $\Sigma(q_1q_2q_3q_4)$.
The edges connecting the circles are labeled by letters
$r_1,...,r_m$ corresponding to the letters of $u$.
The cells tessellating the space between the circles
have labels

\begin{itemize}
\item $q_i^{r_j}=a_jq_i$, $i=2,3,4$, $j=1,...,m$; \  $q_1^{r_j}=q_1$.
\item $ar=ra$, $a\in \{a_1,...,a_m\}$, $r\in  \{r_1,...,r_m\}$
\item $kr=rk$, $k\in\{k_1,...,k_N\}$, $r\in \{r_1,...,r_m\}$.
\end{itemize}

These are exactly the relations of our group $G(v,m)$.

\bigskip

The following two lemmas summarize two main features of the
presentation of the group $G(v,m)$. For any reduced words
$X_1,...X_k$ in the alphabet $a_1^{\pm 1},\dots,a_m^{\pm 1}$ we
define the word $\Lambda (X_1,\dots,X_k)$ to be the word
$u(X_1,...,X_k)$ with letters $q_1,...,q_M$ inserted between
factors $X_{i_j}$ (recall that by the main property of ${\cal
S}(v)$ these are all acceptable words of the $S$-machine ${\cal
S}(v)$.) More precisely $$\Lambda(X_1,\dots,X_k)\equiv
q_1X_{i_1}q_2X_{i_2}...X_{i_{M-1}}q_M,$$ provided
$v(x_1,...,x_k)\equiv x_{i_1}x_{i_2}...x_{i_{M-1}}$. The following
claim is an immediate corollary of the relations (2.1)-(2.5) (and
is evident from Figure 1).

\medskip

\begin{lm} \label{2.1} Assume $i=(j,\ell)$. Then in view of relations
(2.1) - (2.4), the word $r_i^{-1}\Lambda(X_1,\dots,X_k)r_i$ (the word
$r_i\Lambda(X_1,\dots,X_k)r_i^{-1}$) is equal to the word
$\Lambda(X'_1,\dots,X'_k)$, where in the free group
$X'_u = X_u$ for $u\ne \ell$, and $X'_\ell=X_\ell a_j$ ($X'_\ell=X_\ell a_j^{-1}$).
\end{lm}
$\Box$

Now set

$$\Sigma (X_1,\dots,X_k)\equiv k_1\Lambda(X_1,\dots,X_k)\dots
k_N\Lambda(X_1,\dots,X_k).$$

In particular $\Sigma(1,\dots,1)$ is the left-hand side of the
relation (2.6). Since any $k$-tuple $(X_1,\dots,X_k)$ is a result
of iterated multiplications of the components of the $k$-tuple
$(1,\dots,1)$ by letters $a_j^{\pm 1}$, Lemma \ref{2.1} and
relations (2.6) imply

\begin{lm} \label{2.2} For any reduced words $X_1,\dots,X_k$, the word
$\Sigma(X_1,\dots,X_k)$ is conjugate to the
word $\Sigma(1,\dots,1)$ in virtue of relations (2.1) - (2.6).
In particular, $\Sigma(X_1,\dots,X_k)=1$
in the group $G(v,m)$.
\end{lm}
  $\Box$

This lemma is also evident from Figure 1, because Figure 1 is the van Kampen diagram over the presentation of $G(v,m)$ with boundary label $$\Sigma(X_1,...,X_k).$$

\medskip

Now let us define the presentation of the group $H=H(v,m)$ (the
construction is similar to the Aanderaa construction from
\cite{Aanderaa}). Let us add new letters $\rho, d, b_1,\dots,b_m$
to the above presentation of $G(v,m)$, and add the following
relations:

$$\rho^{-1}k_1\rho = k_1 d^{-1},\; \rho^{-1}k_2\rho =dk_2; \eqno (2.7)$$

$$\rho^{-1}k_j\rho = k_j,\; 3\le j\le N; \eqno (2.8)$$

$$\rho^{-1}q_j\rho = q_j,\; 1\le j\le M; \eqno (2.9)$$

$$\rho^{-1}a_j\rho = a_j, \;1\le j\le m; \eqno (2.10)$$

$$d^{-1} a_j d = a_j b_j, \;1\le j\le m; \eqno (2.11)$$

$$d^{-1} q_j d = q_j, \;1\le j\le M; \eqno (2.12)$$

$$ b_j a_\ell = a_\ell b_j, \;1\le j,l\le m; \eqno (2.13)$$

$$ b_j q_\ell = q_\ell b_j, \;1\le j\le m, 1\le l\le M; \eqno (2.14)$$

$$ w(b_1,\dots, b_m) = 1 \eqno (2.15)$$
for any cyclically reduced word
$w=w(b_1,\dots,b_m)$ which is equal to 1 in the relatively free group $F_m({\cal V})$
with the basis $(b_1,\dots,b_m)$.

It is easy to see that these relations do not depend on the
structure of the word $v$ (only on the length of $v$). All these
relations can be put in one picture, the following Figure 2. These
relations together with the relations of $G(v,m)$ form the
presentation of $H(v,m)$.

\unitlength=0.80mm
\special{em:linewidth 0.4pt}
\linethickness{0.4pt}
\begin{picture}(134.66,118.67)
\put(79.50,28.83){\oval(63.00,46.33)[]}
\put(89.50,50.17){\oval(90.33,97.00)[lb]}
\put(98.83,50.83){\oval(32.33,98.33)[rb]}
\put(44.33,47.00){\line(0,1){7.33}}
\put(53.33,5.67){\line(-2,-3){2.67}}
\put(48.00,11.00){\line(-6,-5){3.67}}
\put(48.00,46.67){\line(-1,2){3.67}}
\put(105.33,5.67){\line(1,-1){4.00}}
\put(111.00,11.67){\line(1,-1){4.00}}
\put(80.67,1.67){\line(1,0){23.67}}
\bezier{492}(53.33,63.67)(82.00,118.67)(108.33,63.67)
\bezier{180}(108.67,63.00)(94.67,84.00)(80.00,69.67)
\bezier{136}(108.67,63.00)(91.00,58.34)(80.00,69.67)
\put(78.67,94.67){\makebox(0,0)[cc]{$q_1uq_2uq_3uq_4$}}
\put(90.67,70.67){\makebox(0,0)[cc]{$u_b^3$}}
\put(42.33,60.00){\makebox(0,0)[cc]{$k_1$}}
\put(115.00,61.67){\makebox(0,0)[cc]{$k_2$}}
\put(115.33,1.67){\makebox(0,0)[cc]{$k_3$}}
\put(76.00,48.33){\makebox(0,0)[cc]{$q_1uq_2uq_3uq_4$}}
\put(55.33,54.00){\makebox(0,0)[cc]{$\rho$}}
\put(76.67,25.67){\makebox(0,0)[cc]{Disc}}
\put(53.50,54.50){\oval(18.00,18.50)[lt]}
\put(107.83,51.00){\oval(14.33,24.00)[rt]}
\put(53.33,52.00){\line(0,1){11.67}}
\put(105.33,52.00){\line(1,4){2.67}}
\put(111.00,46.33){\line(1,3){3.67}}
\put(53.33,53.67){\vector(0,1){1.00}}
\put(53.33,58.00){\circle*{0.67}}
\put(53.33,60.33){\vector(0,1){1.67}}
\put(55.00,60.67){\makebox(0,0)[cc]{$d$}}
\put(106.67,57.67){\circle*{0.67}}
\put(105.87,54.33){\vector(1,4){0.67}}
\put(107.18,59.67){\vector(1,4){0.67}}
\put(53.33,58.00){\line(1,0){53.33}}
\put(104.00,55.00){\makebox(0,0)[cc]{$\rho$}}
\put(105.67,60.33){\makebox(0,0)[cc]{$d$}}
\end{picture}

\begin{center}
\nopagebreak[4]
Fig. 2.

\end{center}
\addtocounter{ppp}{1}

This is an annular diagram over the presentation of $H(v,m)$ (as
above we assume for simplicity that $v=x^3$).  It is obtained in
the following way. Take the disc $\Delta$ on Figure 1. Since
$\rho$ commutes with all generators of $G(v,m)$ except $k_1$ and
$k_2$, and $k_1^\rho=k_1d^{-1}, k_2^\rho=dk_2$, we can form an
annulus of $\rho$-cells with the inner boundary labeled by the
same word as the boundary of $\Delta$, and the outer boundary
labeled by the same word with $d^{-1}$ inserted next to the right
of $k_1$ and $d$ inserted next to the left of $k_2$. Glue in the
disc $\Delta$ inside this annulus. We obtain the part of the
diagram on Figure 2 formed by the disc and the $\rho$-annulus
enveloping the disc. Let us call this part $\Delta_1$. Now $d$
commutes with all $q$'s, and we have that $a_i^d=a_ib_i$. Also
take into account that $b_i$ commutes with all the $a$'s and
$q$'s. This implies that if $U=q_1uq_2uq_3uq_4$ is the word
written between $k_1$ and $k_2$ on the boundary of $\Delta$ (read
clockwise) then
$$U^d=q_1uu_bq_2uu_bq_3uu_bq_4=q_1uq_2uq_3uq_4u_b^3=Uu_b^3$$ (all
equalities hold modulo the presentation of $H(v,m)$). Here $u_b$
is the word $u$ rewritten in the alphabet $\{b_1,...,b_m\}$. The
corresponding diagram over $H(v,m)$ can be attached to $\Delta_1$
along the arc labeled by $d^{-1}Wd$. Let the resulting diagram be
denoted by $\Delta_2$.  Now we can get the diagram on Figure 2 by
identifying the ends of the arc labeled by $u_b^3$ on the boundary
of $\Delta_2$.

Notice that the outer boundary of the diagram on Figure 2 is
labeled by the same word as the boundary of the disc $\Delta$.
Thus all the relations $u_b^3=1$ (relations (2.15) above) follow
from the other relations from the presentation of $H(x^3,m)$. The
next statement generalizes these observations to an arbitrary $v$.

\medskip

\begin{lm} \label{2.3} Let $X_1,\dots, X_k$ be words in $a_1^{\pm 1},\dots,a_m^{\pm 1}$, and assume that $Y_1,\ldots,Y_k$ are their copies in $b_1^{\pm 1}, \dots, b_m^{\pm 1}$ (obtained by replacing every $a_j$ by $b_j$ in $X_1,\ldots, X_k$). Then the relation
$$d^{-1}\Lambda(X_1,\dots,X_k)d=\Lambda(X_1,\dots,X_k)v(Y_1,\dots,Y_k)$$
follows from relations (2.11) - (2.14); relations (2.15) follows from
relations (2.1) - (2.14). In particular, $H$ is a finitely presented group.
\end{lm}

\medskip

$\Box$ The first claim is an immediate corollary of the relations
(2.11) -- (2.14) and the definition of the word
$\Lambda\equiv\Lambda(X_1,\dots ,X_k)$. As for the second
statement, it suffices to prove it for all words of the type
$w(b_1,\dots,b_m)\equiv v(Y_1,\dots,Y_k)$ only. In order to do
that, we will apply Lemma \ref{2.2}, relations (2.7) -- (2.10),
and the first claim of Lemma 2.3:

$$1 = \rho^{-1}\Sigma(X_1,\dots,X_k)\rho =
k_1 d^{-1}\Lambda dk_2\Lambda \dots k_N \Lambda= $$
$$k_1 \Lambda v(Y_1,\dots,Y_k)k_2\Lambda\dots k_N\Lambda.$$
By Lemma \ref{2.2} the last product remains being equal to 1 after erasing
the factor $v(Y_1,\dots,Y_k)$. Hence this factor vanishes itself.  $\Box$

\bigskip

\section{Bands and annuli.}

\medskip

Consider a simply connected van Kampen diagram $\Delta$ over the
presentation of the group $H$ (see \cite{LS} or \cite{Ol89}; we
assume that any edge of any van Kampen diagram is labeled by one
letter, as in \cite{Ol89}). If a face $\Pi$ of $\Delta$
corresponds to a relation containing letters $x$ and $y$ then
$\Pi$ is said to be a $(x,y)$-cell. Thus we can talk about
$(\rho,a)$-cells, $(a_j,\rho)$-cell, $(r, q)$-cells, etc.
Similarly if the relation contains letter $x$ then we shall call
the corresponding cell an $x$-cell.

The boundary of a $\rho$-cell $\Pi$ has exactly two $\rho$-edges
labeled by $\rho^{\pm 1}$. These labels are inverses of each other
when one reads the boundary label of $\Pi$. This gives us an
opportunity to construct ``bands" of several $\rho$-cells.

A $\rho$-band of length 0 has no faces and consists of one
$\rho$-edge. A $\rho$-band of length 1 is just a single
$\rho$-cell.  Assume by induction, that we have a $\rho$-band $T'
= [\Pi_1,\dots,\Pi_{s-1}]$ of length $s-1$ constructed of $s-1$
distinct $\rho$-cells $\Pi_1,\dots,\Pi_{s-1}$, and the boundary of
$T'$ has a $\rho$-edge $e$, which is a common edge of the boundary
$\partial(\Pi_{s-1})$ and the boundary of a $\rho$-cell $\Pi_s$
which is distinct from $\Pi_1,\dots,\Pi_{s-1}$. Then we are able
to construct a $\rho$-band $T= [\Pi_1,\dots,\Pi_s]$ of length $s$
whose boundary $\partial T$ is the union of the  the boundaries
$\partial T'$ and $\partial\Pi_s$ minus the edge $e$. A
$\rho$-band $T$ is {\it maximal} if is not contained in a
$\rho$-band of a greater length.

Thus, the boundary of a $\rho$-band $T$ has the form $e_1pe_2q$, where
$e_1$ and $e_2$ are $\rho$-edges (we call them {\it ends} of $T$), and
the paths $p, q$ ({\it sides} of the band $T$) consist of $a$- and
$q$-edges, but contain no $\rho$-edges.

If the ends of $T$ coincide, one may identify them and the annular
subdiagram $T$ is called a {\it $\rho$-annulus }. For example the
diagram on Figure 2 contains a $\rho$-annulus enveloping the disc.

Similarly we define $d$-bands (and annuli) that by definition can be
constructed of $(d,a)$-, $(d,t)$-, and $(d,q)$-cells.

$b$-bands can be constructed of $(b,a)$- and $(b,q)$-cells.

$r$-bands are constructed of $(r,q)$-, $(r,a)$-, and $(r,k)$-cells.

$q$-bands are constructed of $(q,r)$-, $(q,\rho)$-, $(q,d)$-, and $(q,b)$-cells.

$a$-bands are created of $(a,r)$-, $(a,\rho)$-, $(a,d)$-, and $(a,b)$-cells.

Notice that $(\rho,d)$-cells of type (2.7) cannot be included in a
$d$-band but they can be {\it terminal} for $d$-bands, i.e. a maximal
$d$-band can end only on the contour of a $(\rho,k)$-cell or on
the contour of the diagram $\Delta$.  Similarly, $(r,q)$-cells
are terminal for $a$-bands, {\it hubs}, corresponding to relation (2.6),
are terminal for $k$- and $q$-bands, $(d,a)$-cells are
terminal for $b$-bands. Also a $b$-band can terminate on the
contour of a $G_b$-cell (by definition, a $G_b$-{\it cell} corresponds
to a relation (2.15)).

Now consider an $r$-band $T=[\pi_0,\pi_1,\dots,\pi_\ell,\pi_{\ell+1}]$ and
a $q$-band \\$T'=[\pi_0,\gamma_1,\dots,\gamma_s,\pi_{s+1}]$ which have
no common faces except for $\pi_0$ and $\pi_{\ell+1}$, and with all ends
of $T$ and $T'$ lying on the outer boundary of the annulus $S$ formed by $T$
and $T'$. Then this annulus is called an $(r,q)-${\it annulus}.
It consists of the $r$-{\it part} $T$ and the $q$-{\it part} $T'$. The faces
$\pi_0$ and $\pi_{\ell+1}$ are its {\it corner} cells.

The definitions of $(\rho,a)$-, $(q,b)$-annuli, etc. are quite similar.

A diagram $\Delta$ is called {\it minimal} in this section
if there exists no other diagram $\Delta'$ such that (1) $\Delta'$ has
the same boundary label as $\Delta$, and (2) the number of faces of each of
the types (2.1) - (2.15)
in $\Delta'$ does not exceed the similar number for $\Delta$, (3) the
total number of faces in $\Delta'$ is smaller than the number of faces
in $\Delta$. In the next two sections we shall a stronger definition of minimality.

The main lemma of this section claims that there are no annuli of
various kinds in minimal diagrams without hubs.

\medskip


\begin{lm}\label{3.1} Let $\Delta$ be a minimal diagram over $H$
containing no hubs. Then $\Delta$ has no

(1) $\rho$-annuli,

(2) $r$-annuli,

(3) $(r,q)$-annuli,

(4) $q$-annuli,

(5) $(r,k)$-annuli,

(6) $k$-annuli,

(7) $(\rho,k)$-annuli,

(8) $(a,b)$-annuli,

(9) $(d,a)$-annuli,

(10) $(\rho,a)$-annuli,

(11) $(r,a)$-annuli,

(12) $d$-annuli,

(13) $b$-annuli,

(14) $a$-annuli,

(15) $(\rho,q)$-annuli,

(16) $(d,q)$-annuli,

(17) $(q,b)$-annuli
\end{lm}
\medskip

$\Box$ To prove statements (1) - (17) we use a simultaneous induction
on the number of faces in the minimal subdiagram $\Delta_S$
containing a conjectural counterexample, i.e. an annulus $S$. This
means that we may assume that $\Delta$ has no annulus $S'$ of types
(1) - (17) such that the subdiagram $\Delta_{S'}$ has fewer faces
than $\Delta_S$.

(1) Let $S$ be a $\rho$-annulus. Assume that $S$ has a
$(\rho,k_j)$-cell. Then this cell belongs to a $k_j$-band $T$,
which must intersect $S$ at least twice, because by the lemma
condition $\Delta$ contains no hubs (terminal cells for
$k$-bands). In such a case $T$ and a subband of $S$ form a smaller
$(\rho,k)$-annulus $S'$ than $S$, contrary to claim (7) of the
lemma. The only case when $S'$ is not smaller than $S$ is when $S$
consists just of two $(\rho,k_j)$-cells with a common $\rho$-edge
and with ``mirror" labels. Such a pair of mirror faces is
impossible in a minimal diagram. For more details on the cell
cancellation, see \cite{LS}.

Therefore $S$ contains only $(\rho,a)$- and $(\rho,q)$-cells.
Consequently the outer and the inner boundaries of $S$ have identical
labels. This makes it possible to delete the interior of $S$
and then identify the sides of $S$. Such a surgery does
not change the boundary label of $\Delta$, contrary to the
minimality of $\Delta$.

(3) Assume that $S$ is a $(r,q)$-annulus. Let $T$ be the $q$-part
and $T'$ be the $r$-part. $S$ has no $k$-cells, because otherwise
a smaller $(r,k)$-annulus appears, contradicting claim (5).
Analogously, by (3) (for smaller
annuli) there are no non-corner $(r,q)$-cells in $S$. The same
argument shows that there are no other $(r,q)$-cells in the
subdiagram $\Delta $.

The corner $(r,q)$-cells are included in the same $r$- and $q$-bands.
This implies that they correspond to the same relation and
simultaneously have or have no $a_j$-edges (for the same $j$) on
the inner border of $S$ with opposite directions. Only these
corner cells can be terminal for $a$-bands crossing $S$.

Therefore if there exist non-corner cells in $S$, then by (11) there
exist exactly 2 such cells, they must be neighbors in the
$r$-part of $S$, and must have ``mirror" labels, contrary to the
minimality of $\Delta $.

Hence the $r$-part of $S$ has no non-corner cells. Then the standard
cancelation argument can be applied to the corner cell. This
contradicts the minimality assumption again.

(2), (4) - (17) The proofs of all these statements are similar to
the two proofs of (1) and (3) given above. The reader could
examine them as an exercise or read similar explanation for Lemma
6.1 \cite{Ol1}  (claims (1) - (20)) or Section 7 of \cite{SBR}.
$\Box$

 \bigskip

\section{Hubs and spokes}

\medskip

A {\it spoke} is a maximal $k$-band having an end on a hub (or on
a disc in the next section). Obviously, another end lies on
a hub too, or on the boundary $\partial\Delta$.

It will be convenient to restrict the notion of a minimal diagram
used in previous section as follows.

A {\it type} $\tau = \tau(\Delta)$ of a diagram $\Delta$ over $H$
is the 4-tuple $\tau = (\tau_1,\tau_2,\tau_3,\tau_4)$ where
$\tau_1$ is the number of hubs in $\Delta$ (or, in the next
section, the number of discs), $\tau_2$ is the number of
$(\rho,k)$- and $(r,k)$-cells, $\tau_3$ is the number of all other
faces except $(b,q)$-, and $(b,a)$- cells, $\tau_4$ is the number
of $(b,t)$-, $(b,q)$-, and $(b,a)$-cells in $\Delta$. Set $\tau <
\tau'$ if $\tau_1 <\tau'_1$, or $\tau_1 = \tau'_1$, but $\tau_2
<\tau'_2$, and so on. Further a diagram is said to be {\it
minimal} if it has the minimal type among all diagrams with the
same boundary label. Clearly, this notion of minimality is
stronger than that in Section 3, that is a diagram which is
minimal under this definition is also minimal under the definition
in the previous section.

\medskip

Since our construction of the group $G(v,m)$ is essentially the
same as the construction of the group $G_N({\cal S})$ of
\cite{SBR}, the next lemma follows from Lemma 11.1 of \cite{SBR}.
Nevertheless we present a direct proof here.

\begin{lm}\label{4.1} Let $\Delta$ be a minimal diagram over $G(v,m)$
containing two hubs $\Pi_1$ and $\Pi_2$. Assume that the hubs have
two common consecutive spokes $T_1$ and $T_2$ such that the subdiagram $\Delta_0$
bounded by these spokes and the hubs contains no hubs. Then
$\Delta$ is not minimal diagram: by passing to another diagram with
the same boundary label one can
decrease the number of hubs by 2.
\end{lm}
\medskip

$\Box$ Let $\Delta_1$ be the subdiagram constructed of $\Delta_0$
and the spokes $T_1$, $T_2$. Without loss of generality we may
assume that $\Delta_1$ is a minimal diagram. By Lemma \ref{3.1}
(5) for $\Delta_1$, every $r$-band crossing $T_1$, must cross
$T_2$ as well. Therefore each of the 4 sides of the bands $T_1$,
$T_2$ must have the same label $V = V(r_1,r_2,\dots)$. Thus the
boundary label of $\Delta_0$ is the commutator $[V,\Lambda_1]$ for
the word $\Lambda_1\equiv \Lambda(1,\dots,1)\equiv q_1tq_2t\dots
q_M$, i.e. $V$ commutes with $\Lambda_1$ modulo all the defining
relations of $G(v,m)$ excluding the hub. Being a word in the
alphabet $r_1^{\pm 1}, r_2^{\pm 1},\dots$, $V$ commutes with
letters $k_j$ (see relations (2.5)). Consequently, $V$ commutes
(modulo (2.1) - (2.5)) with the word $W$, which is a cyclic
permutation of the left-hand side of the hub relation (2.6)
written on $\Pi_1$, $\Pi_2$ in opposite directions starting with
the ends of the band $T_1$.

Thus, if we cut out the hubs $\Pi_1$, $\Pi_2$ from $\Delta$ and then make a cut
along the band $T_1$ border, we get a hole labeled by the word $[V, W]$
equal to 1 by (2.1) -- (2.5). Now by the van Kampen lemma we are able to
insert a diagram of a type $(0,\star,\star,\star)$ in this hole, reducing the number of
hubs in $\Delta$ by 2.  $\Box$

\medskip

\begin{lm}\label{4.2} Let $\Delta$ be a diagram over $H(v,m)$ containing
two hubs $\Pi_1$, $\Pi_2$ with common consecutive $k_j$- and $k_s$-spokes
$T_1$ and $T_2$ where $\{j,s\}\ne \{1,2\}$. Assume that the subdiagram
$\Delta_0$ bounded by these spokes and the hubs contains no hubs. Then
the diagram $\Delta$ is not minimal: By passing to another diagram with
the same boundary label, one can decrease the number of hubs by 2.
\end{lm}
\medskip

$\Box$ The proof is similar to the proof of Lemma \ref{4.1}, but
now the word $V$ may contain $\rho^{\pm 1}, r_1, r_2,\dots $. So,
the difference is that $V$ may not commute with the letters $k_1$,
$k_2$ occurring in $W$. However $k_1, k_2$ occur just in the
subword  $k_1\Lambda_1k_2$ of the left-hand side in (2.7), which
is also a subword of $W$ in view of the condition $\{k_j,k_s\}\ne
\{1,2\}$. Therefore it suffices to check that $V$ commutes with
$k_\ell$ for $\ell\ne 1,2$, with $\Lambda_1$, and with
$k_1\Lambda_1k_2$ modulo relations (2.1) - (2.5), (2.7) - (2.15).
The first property follows from (2.5) and (2.9). The second one
can be explained exactly as in Lemma \ref{4.1}.

To prove the third commutativity, notice first of all that by
relations (2.15) and Lemma \ref{2.3} $d$ commutes with
$\Lambda(X_1,\dots,X_k)$ for any words $X_1,\dots,X_k$ in the
alphabet $\{a_1^{\pm 1},\dots,a_m^{\pm 1}\}$. Then $\rho$ commutes
with both $\Lambda(X_1,\dots,X_k)$ (see (2.9), (2.10)) and
$k_1\Lambda(X_1,\dots,X_k)k_2$ since by (2.7)
$\rho^{-1}k_1=k_1d^{-1}\rho^{-1}$ and $k_2\rho=\rho dk_2$.

Also recall that by (2.5) and Lemma \ref{2.2}
$$r_i^{\pm 1}k_1\Lambda(X_1,\dots,X_k)k_2 r_i^{\mp 1}=
k_1 r_i^{\pm 1}\Lambda(X_1,\dots,X_k)r_i^{\mp 1}k_2 =
k_1\Lambda(X'_1,\dots,X'_k)k_2$$
for some $X'_1,\dots,X'_k$.

Therefore permuting the words $V^{-1}$ and $V$
(these are word in in $\rho^{\pm 1}, r_1^{\pm1}, r_2^{\pm 1},\dots $)
letter-by-letter with
$k_1$ and $k_2$,
we get
$$V^{-1}k_1\Lambda_1 k_2 V =k_1 V^{-1}\Lambda_1 V k_2.$$
The right-hand side is equal to $k_1\Lambda_1 k_2$, as was mentioned earlier. $\Box$

\medskip

With any minimal diagram $\Delta$ over $H(v,m)$ we associate the following
graph $\Gamma=\Gamma_{\Delta}$. One vertex of $\Gamma$ ({\it exterior}
vertex) is taken outside
$\Delta$ on the plane. Every {\it interior} vertex is chosen inside
a hub. The edges between the vertices are drawn along the "medians"
of the spokes. The {\it exterior} edges are incident to the exterior
vertex. The other edges are {\it interior}. Finally to complete the
definition of $\Gamma$, we erase the interior $k_1$-spoke for every
pair of hubs connected by both $k_1$- and $k_2$-spokes.

By Lemma \ref{4.2} there are no bigons formed by interior edges of
$\Gamma$. Also there are no loops in $\Gamma$ because every letter
$k_j$ occurs once in the left-hand side of (2.7). Therefore in
standard way the Euler formula implies that $\Gamma$ has many
exterior edges. (The restriction $N-1\ge 6$ would be enough here;
the stronger condition $N-1\ge 28$ will be useful for Sections
5-7.) The following statement (see  Lemma 2.13 in \cite{BORS}, or
Lemmas 3.2, 3.3 in \cite{Ol1}, or Lemma 11.5 in \cite{SBR}) will
be sufficient for our purpose.

\medskip

\begin{lm}\label{4.3} Let $\Delta$ be a minimal diagram containing
at least one hub. Then there exists a hub $\Pi$ in $\Delta$ such
that at least $N-4$ consecutive spokes starting on $\Pi$, have
their ends on the boundary $\partial\Delta$, and moreover there
are no other hubs between the spokes of this set. The number of
$k$-edges in $\partial\Delta$ is at least 3 times greater than the
number of hubs in $\Delta$.\end{lm} $\Box$

\medskip

\begin{lm}\label{4.4}The natural homomorphism of the group $F_m({\cal V})$
onto $H(v,m)$ (well defined in view of relations (2.15)) is injective.
\end{lm}

\medskip

$\Box$  Assume that a word $w\equiv w(b_1,\dots,b_2)$ vanishes
under the homomorphism. Then there exists a minimal van Kampen
diagram $\Delta$ with the boundary label $w$. By Lemma \ref{4.3}
there are no hubs in $\Delta$ because no letter $k_j$ occurs in
$w$. Then by Lemma \ref{3.1} (1) $\Delta$ has no $\rho$-annuli,
and consequently, it has no $\rho$-cells at all. Quite similarly,
Lemma \ref{3.1} allows us to exclude $r$-, $q$-, $k$-, $d$-, and
$a$-cells from $\Delta$ consequently. For example, there are no
$d$-cells because maximal $d$-bands could terminate on
$\rho$-cells only. Thus, $\Delta$ has $G_b$-cells only, that
correspond to relations (2.15). Hence $w=1$ in $G$ as desired.
$\Box$

\bigskip

\section{The band structure of disc-based diagrams}

\medskip

It will be convenient to extend the list of relations of the group
$H(v,m)$ by adding the relations $\Sigma(X_1,\dots,X_k)=1$ for all $k$-tuples
of reduced words $X_1,\dots,X_k$ in $a_1^{\pm 1},\dots,a_m^{\pm 1}$.
Such an enlargement does not change $H(v,m)$ by Lemma \ref{2.2}. Any face
of a diagram over this presentation of $H(v,m)$, that corresponds to some
relation $\Sigma(X_1,\dots,X_k)=1$, will be called  a {\it disc}.

The notion of a minimal diagram will be further restricted by
replacing discs for hubs in the definition of a minimal diagram
from the previous section. With any minimal diagram $\Delta$ we
associate a graph $\Gamma(\Delta)$. Its definition repeats the
definition of the graph $\Gamma_{\Delta}$ given in Section 4 where
discs replace hubs. Notice that by Lemma \ref{2.2} every disc can
be replaced by a hub (which is a disc too) and a number of faces
of smaller ranks  (but the resulted diagram may be not minimal).
The possibility of such a replacement and Lemma \ref{4.2} show
that the graph $\Gamma(\Delta)$ has no bigons as well. Therefore
the statement of Lemma \ref{4.3} is also true for
$\Gamma(\Delta)$.

   Let $S$ be a $\rho$- or $r$-band that consequently intersects
at $(\rho,k)$- or $(r,k)$-cells a series of consequent spokes
$T_1,\dots,T_\ell$ starting on a disc $D$. We say that the band $S$
{\it envelopes} disc $D$ if $\ell>N/2$, and there are no other discs in
the sectors formed by $S$, $T_j$ and $T_{j+1}$ for $j=1,\dots,\ell-1$.

\medskip

\begin{lm} \label{5.1} A minimal diagram $\Delta$ over $H(v,m)$ has no
bands which envelope discs.
\end{lm}

\medskip

$\Box$ The proof is completely similar to the proofs of Lemma 8.4
\cite{Ol1} or Lemma 2.17 in \cite{BORS}. Therefore we give just a brief
explanation below referring for details to \cite{Ol1} or \cite{BORS}.

Arguing by contradiction, one can choose the closest to $D$ band $S$
that envelopes it. Then the intersection cells $\pi_1,\dots,\pi_\ell$ of $S$
and $T_1,\dots,T_\ell$ have common $k$-edges with $D$. It suffices to prove
that the diagram $\Delta_0$ consisting of $D$ and $\pi_1,\dots,\pi_\ell$
and considered  separately from $\Delta$, is not minimal.

For this purpose we attach auxiliary $(r_s,k)$-cells $\pi_{\ell+1},\dots,\pi_N$
to $\Delta_0$ along the $N-\ell$ free $k$-edges of the disc $D$ (if $S$ is a
$r_s$ band) so that all the cells $\pi_1,\dots,\pi_N$ would be attached to $D$
uniformly, i.e. their $r$-edges would be directed all 'to' or all 'from' $D$.
Then adding several faces of smaller ranks we can get a diagram $\Delta_1$
with a label $\Sigma(X'_1,\dots,X'_k)$ by Lemma \ref{2.2}.

Therefore, conversely, one can construct a diagram $\Delta_2$, with the same
boundary label as $\Delta_0$, consisting of a disc $D'$ labeled by $\Sigma(X'_1,\dots,X'_k)$,
mirror copies of $(k,r)$-cells $\pi_{\ell+1},\dots,\pi_N$, and faces of smaller
ranks. But this contradicts to the minimality of $\Delta_0$ since $N-\ell<\ell$.

If $S$ is a $\rho$-band, the proof is similar, but the boundary label of the
disc $D'$ coincides with the label of $D$, since $\rho$ commutes with
$\Sigma(X_1,\dots,X_k)$ as was explained in Lemma \ref{4.2}.  $\Box$

\medskip

\begin{lm}\label{5.2} A minimal diagram $\Delta $ has no annuli of types
(1) - (17) from Lemma 3.1 (even if hubs or discs occur in $\Delta$).\end{lm}

\medskip

$\Box$ This is similar to Lemmas in Section 4 of \cite{BORS}. Let
us show, for example, that statement (3) of this lemma can be
deduced from statement (3) of Lemma \ref{3.1}.

Let $\Delta_S$ be a minimal subdiagram of $\Delta$, containing a
$(r,q)$-annulus $S$. By Lemma 3.1(3) it contains a disc. By Lemma 4.3
there exists a disc $D$ in $\Delta_S$ such that its consecutive spokes
$T_1,\dots,T_{N-4}$ intersect the $r$-part $R$ of $S$ (since the $q$-part
has no $k$-cells at all). If there are other discs in $\Delta_S$
and the $q$-part $Q$ of $S$ occurs between some $T_j$ and $T_{j+1}$,
then again as in Lemma \ref{4.3} (but for the graph obtained by erasing
the vertex in $D$ and the edges incident to it), we get another
disc $D'$ and spokes $T'_1,\dots,T'_{N-5}$ starting on it, such
that $R$ intersects them consecutively , and there are neither hubs
nor $Q$ between the spokes. But this contradicts Lemma \ref{5.1} because
$N-5\ge N/2 $. $\Box $

\medskip

\begin{lm} \label{5.3} Any two distinct maximal bands $T$ and $T'$
have at most one common face in a minimal diagram over $H(v,m)$.
\end{lm}
\medskip

$\Box$ Basically the statement follows from Lemma \ref{5.2}.
However we have to remember that a-priori, a multiple intersection
of two bands $T$ and $T'$  does not imply that they form even one
annulus, because one or both ends of the band can be inside the
``annulus". Figure 3 shows a spiral multiple intersection.

\begin{center}
\unitlength=1.50mm
\linethickness{0.4pt}
\begin{picture}(97.67,32.89)
\put(11.67,10.67){\framebox(24.00,3.56)[cc]{}}
\bezier{88}(35.67,14.22)(38.78,24.67)(42.33,14.22)
\bezier{200}(42.33,14.22)(45.22,-2.67)(15.22,10.67)
\bezier{168}(32.33,14.22)(39.22,32.89)(46.11,12.22)
\bezier{240}(46.11,12.22)(48.33,-6.89)(11.67,10.67)
\put(15.22,10.67){\line(0,1){3.56}}
\put(32.33,14.22){\line(0,-1){3.56}}
\put(13.44,12.22){\makebox(0,0)[cc]{$\pi$}}
\put(34.11,12.45){\makebox(0,0)[cc]{$\pi'$}}
\put(24.11,12.45){\makebox(0,0)[cc]{$\ldots$}}
\put(36.78,17.33){\line(-1,1){2.40}}
\put(43.44,16.00){\circle*{0.00}}
\put(44.33,12.67){\circle*{0.00}}
\put(43.89,9.11){\circle*{0.00}}
\put(18.56,14.22){\line(0,-1){3.56}}
\put(68.78,10.67){\framebox(20.44,3.56)[cc]{}}
\put(72.33,10.67){\line(0,1){3.56}}
\put(72.33,14.22){\line(0,0){0.00}}
\put(75.67,14.22){\line(0,-1){3.56}}
\put(85.00,14.22){\line(0,-1){3.56}}
\put(70.56,12.44){\makebox(0,0)[cc]{$\pi$}}
\put(87.22,12.22){\makebox(0,0)[cc]{$\pi'$}}
\bezier{80}(89.22,14.22)(92.33,22.00)(93.00,10.67)
\bezier{80}(93.00,10.67)(94.11,6.44)(78.78,6.44)
\bezier{100}(84.33,6.67)(65.00,4.89)(65.22,10.67)
\bezier{72}(65.22,10.89)(65.22,20.89)(68.78,14.22)
\bezier{168}(85.00,14.22)(95.00,32.00)(96.11,10.67)
\bezier{108}(96.11,10.67)(97.67,2.89)(78.78,2.89)
\bezier{128}(83.89,3.11)(62.11,0.89)(62.33,10.67)
\bezier{156}(62.33,10.67)(63.89,30.89)(72.33,14.22)
\put(90.11,16.22){\line(-1,1){2.67}}
\put(66.33,16.89){\line(1,1){2.89}}
\put(80.56,12.22){\makebox(0,0)[cc]{$\ldots$}}
\put(81.00,4.67){\makebox(0,0)[cc]{$\dots$}}
\end{picture}
\end{center}

\begin{center}
\nopagebreak[4]
Fig. 3.
\end{center}

In fact such spirals cannot occur, and the reader can find details
in Lemmas 5.1 and 5.8 from \cite{BORS}. Here we just explain the
idea for the particular case when $T$ is a $d$-band and $T'$ is an
$a$-band. In this case a terminal cell $\Pi$ for $T$ must be
inside the subdiagram $\Delta_S$ bounded by $S$. $\Pi$ is a
$(\rho,k)$-cell (see relations (2.7)). Therefore by Lemma
\ref{4.3} and Lemma \ref{3.1} (6) the boundary $\partial\Delta_S$
must be crossed by at least one $k$-band. But both $d$-band $T$
and $a$-band $T'$ have no $k$-cells at all, a contradiction.
$\Box$

\medskip

The following statement is similar to Lemma 4.36 in \cite{BORS}.

\medskip

\begin{lm}\label{5.4} There is no $b$-band $T$ in a minimal diagram $\Delta$
such that both ends of $T$ belong to $G_b$-cells.
\end{lm}

\medskip

$\Box$ First assume that both ends of $T$ belong to the boundary
of one $G_b$-cell $\Pi$. Then the boundary of $T$ and a subpath of
$\partial\Pi$ form the boundary of a subdiagram $\Delta_0$ with
the boundary label $UV$ where $U$ is a word in $b_1^{\pm 1},\dots,b_m^{\pm 1}$,
and $V$ is a word in $a_1^{\pm 1},\dots,a_m^{\pm 1}$,
$q_1^{\pm 1},\dots,q_M^{\pm 1}$. By lemmas \ref{4.3} and \ref{3.1}
$\Delta_0$
has no discs and $k$-cells. By lemma \ref{3.1} the band $T$ has
no $a$- and $q$-cells because otherwise we would get
$(a,b)$- or $(q,b)$-bands.
Thus the band $T$ has length 0. Hence there is a loop in $\partial\Pi$
whose label must be equal to 1 in the group $G$ by Lemma \ref{4.4}. Then
one can replace the interior of this loop and $\Pi$ by a single
$G_b$-cell, contrary to the minimality of $\Delta$.

Assume now that the ends of $T$ belong to distinct $G_b$-cells
$\Pi_1$ and $\Pi_2$. Recall that the boundary label of $T$ commutes
with any word in $b_1^{\pm 1},\dots,b_m^{\pm 1}$ by relations (2.13) and (2.14).
Therefore the diagram $\Delta$ is not minimal: the subdiagram consisting
of the two $G_b$-cells connected by the border of $T$, can be replaced
by a diagram consisting of one $G_b$-cell and several cells of smaller
ranks that correspond to relations (2.13) and (2.14).  $\Box$

\bigskip

\section{Upper bounds for the number of cells in minimal diagrams}

\medskip

Recall that by Lemma \ref{2.3} the group $H(v,m)$ has a finite presentation.

\medskip

\begin{lm}\label{6.1} For any word $w$ of length $n$ in the generators of
the group $H(v,m)$, that is equal to 1 in $H(v,m)$, there is a diagram over
the finite presentation (2.1) - (2.14) with $n^2 O(f(O(n^2))^2$ cells,
where $f$ is the function given in Theorem 1.
\end{lm}

\medskip

$\Box$ First consider a minimal diagram $\Delta$
with the boundary label $w$ over
the disc-based presentation.

It follows from the statement of Lemma \ref{4.3} for discs, that the number of
discs in $\Delta$ is $O(n)$. Therefore by Lemma 5.2 the
number of maximal $\rho$-, $r$-, $q$-, and $k$-bands in $\Delta$
is $O(n)$ because the boundary of any disc has no $\rho$- and $r$-edges,
and has $O(1)$ $q$-, and $k$-edges.

Therefore by Lemma \ref{5.3} the number of $(r,q)$-cells (the
intersections of $r$- and $q$-bands) is $O(n^2)$. Similarly, the
number of $(r,k)$-, $(\rho,k)$-, and $(\rho,q)$-cells is $O(n^2)$.

By Lemma \ref{5.2} the number of maximal $d$-bands in $\Delta$ is
$O(n^2)$ since only $(\rho,k)$-cells can be terminal for
$d$-bands.

A similar argument shows that the number $N_1$ of all maximal
$a$-bands terminating on $(r,q)$-cells or on $\partial\Delta$ is
$O(n^2)$. To obtain similar upper bound for the number of all
$a$-bands, it suffices to explain why the number $N_2$ of maximal
$a$-bands having both ends on discs, is not greater than $N_1$.
Indeed the number of common spokes of two discs is at most 2.
Consequently the number of common $a$-bands for 2 discs is at most
$3/N < 1/9$ of the number of maximal $a$-bands starting on each of the disc.
This implies that $N_2\le N_1$.

The upper bounds for the numbers of maximal $k$-, $q$-, and $a$-bands
show that the sum of perimeters of all discs in $\Delta$ is $O(n^2)$.

Now by Lemma \ref{5.3} we are able to conclude that the number of
$(\rho,a)$-, $(r,a)$-, $(d,q)$-, $(d,t)$-cells is $O(n^3)$, and
the number of $(d,a)$-cells is $O(n^4)$.

Only $(d,a)$-cells can be terminal for $b$-bands. So the above
estimate of the number of $(d,a)$-cells together with Lemma
\ref{5.2} imply that the number of maximal $b$-bands is $O(n^4)$.
This implies, first of all, that the number of $(b,a)$-cells is
$O(n^6)$, and the number of $(b,q)$- and $(b,t)$-cells is
$O(n^4)$. Second, this implies that the sum of perimeters of all
$G_b$-cells is $O(n^4)$. By the definition of the function $f$
every $G_b$-cell with boundary length $\le s$ can be tiled by
cells with boundary labels of the form $v(Y_1,\dots,Y_k)$ and with
$\sum |Y_{ij}|$ at most $f(s)$. By the superadditivity of the
function $f$ (see Proposition 1.1) the set of all $G_b$-cells in
$\Delta$ can be replaced by the set of $G_b$-cells with labels of
the form $v(Y_1,\dots,Y_k)$ and with $\sum |Y_{ij}|$ at most
$f(O(n^4))$.

By Lemma \ref{2.3} every $G_b$-cell $\Pi$ labeled by a word $v(Y_1,\dots,Y_k)$,
can be replaced by a some number $N_{\Pi}$ of cells corresponding to
relations (2.1) - (2.14). The straightforward computation shows that
$N_{\Pi} = O(|Y_1|+\dots+|Y_k|)^2$. Therefore the set of all cells
corresponding to relations (2.1) - (2.14), that replace all $G_b$-cells
of $\Delta$ is $O(f(O(n^4))^2)$.

Finally, every disc which is not a hub, can be replaced by a hub and
several cells corresponding to relations (2.1) - (2.6). We need at most
$O(l^2)$ such cells for every disc of perimeter $\ell$ which can be easily
derived from the proof of Lemma \ref{2.2}.
Therefore all discs in $\Delta$ can be replaced by
$O(n^4)$ cells of types (2.1) - (2.5).

Summing all the upper bounds for the numbers of cells of different
types, we obtain an isoperimetric function of the group $H(v,m)$
that is equal to $O(f(O(n^4))^2)$. To achieve the better upper
bound of $n^2 O(f(n^2)^2)$ claimed in the Theorem, we prove that
the perimeter of every $G_b$-cell in $\Delta$ is in fact bounded
by $O(n^2)$. The explanation is completely similar to the proof of
Lemma 5.15 in \cite{BORS}. It is based on the fact that the number
of the maximal $b$-bands starting on the same $G_b$-cell and
terminating on $a$-cells lying in the same $a$-band, cannot exceed
3 (see Lemma 4.34 in \cite{BORS}). We leave details to the reader.
$\Box$

\medskip

{\large Proof of Theorem 1.} $\Box$ The first statement of the theorem
follows from Lemmas 2.3, 4.4 and 6.1. To prove that the constructed
embedding of $G$ into $H(v,m)$ is undistorted, we need a longer and more complicated
argument. The reasoning essentially coincides with that in Sections 10-12
of \cite{Ol1} or in Section 7 in \cite{BORS}.
Therefore we will not repeat it here
referring the reader to \cite{Ol1} and \cite{BORS}. $\Box$

\bigskip

\section{The embeddings of the groups $BS_{k,1}$}

\medskip

It is more convenient to change the names of generators of
$BS_{k,1}$ to $b_1, b_2$. So  $BS_{k,1}=\langle b_1, b_2 \ | \
b_2^{b_1}=b_2^k\rangle$. Obviously the words $W_n=W_n(b_1,b_2)
\equiv (b_1^nb_2b_1^{-n})b_2(b_1^nb_2^{-1}b_1^{-n})b_2^{-1}$
represent the identity in the group $B=BS_{k,1}$, and it is known
that their ``areas" exponentially increase depending on $n$ if
$k\ge 2$. To prove Theorem 2 we are going to embed $BS_{k,1}$ into
a finitely presented group $H=H_{k,1}$ such that areas of the
words $W_n$ with respect of defining relations of $H$ have
quadratic growth (a relatively easier task), and then we will
prove a polynomial isoperimetric inequality for all words
vanishing in $H$ (a harder job).

The embedding is similar to the one used for relatively free
groups (see the proof of Theorem 1 above). But in order to show
flexibility of our approach, we modify the second step of the
embedding slightly. The main difference of the construction that
we are about to present and the construction presented above is
the absence of the letter $d$. Instead, we have different letters
in different sectors of the discs.

An auxiliary group $G=G_{k,1}$ is given by generators $$a_1, a_2,
c, q_1, q_2, q_3, q_4, k_1,\dots,k_N$$ where $N\ge 29$ and by
defining relations

$$rq_1r^{-1}=q_1a_1, rq_2r^{-1}=q_2a_1^{-1}, rq_3r^{-1}=q_3a_1,
rq_4r^{-1}=q_4a_1^{-1}, rq_5r^{-1}=q_5c; \eqno (7.1)$$

$$rk_jr^{-1}=k_j , \;  1\le j \le N; \eqno (7.2)$$

$$ra_1r^{-1}=a_1, \;ra_2r^{-1}=a_2,\; rcr^{-1}=c; \eqno (7.3)$$

$$q_1a_2q_2a_2q_3a_2^{-1}q_4a_2^{-1}k_1q_5\dots k_Nq_5=1. \eqno (7.4)$$

The relation (7.4) is called the {\it hub} relation.

Note that it is easy to draw the disc corresponding to these
relations and write the rules of the corresponding $S$-machine. It
is a good exercise for a reader who wants to learn how to draw van
Kampen diagrams and write programs for $S$-machines.

The words $\Sigma_s=\Sigma (a_1^s, a_2)$ are defined as follows:
$$\Sigma_s\equiv
q_1a_1^sa_2q_2a_1^{-s}a_2q_3a_1^{s}a_2^{-1}q_4a_1^{-s}a_2^{-1}k_1q_5c^s\dots
k_Nq_5 c^s.$$ An obvious analog of Lemma \ref{2.2} says that
$r\Sigma_sr^{-1}=\Sigma_{s+1}$ in view of relations (7.1) --
(7.3), and the deduction takes $O(n^2)$ of application of
relations (7.1) -- (7.3).

The group $H=H_{k,1}$ is defined by adding to the presentation of the
group $G$ new generators $b_1, b_2, \rho$, and by adding new relations

$$\rho a_j\rho^{-1}=a_jb_j,\; j=1,2; \eqno (7.5)$$

$$\rho c\rho^{-1} =c; \eqno (7.6)$$

$$\rho q_j\rho^{-1}=q_j, \;j=1,2,3,4,5; \eqno (7.7)$$

$$\rho k_j\rho^{-1}=k_j, \;1\le j\le N; \eqno (7.8)$$

$$b_j a_i = a_ib_j, \;i,j =1,2; \eqno (7.9)$$

$$b_j q_i =q_i b_j, \;i=1,2,3,4, j=1,2; \eqno (7.10)$$

$$b_1b_2b_1^{-1}=b^k. \eqno (7.11)$$

By definition, a $B$-cell in a diagram over $H$ corresponds to
any cyclically reduced word $w(b_1,b_2)$ such that $w=1$ follows from (7.11).
The following analog of Lemma \ref{2.3} can be verified immediately (and similar to the
proof of Lemma \ref{2.3}).

\medskip

\begin{lm} \label{7.1} For any $n\ge 0$ the relation
$$\Sigma_n\equiv
q_1a_1^na_2q_2a_1^{-n}a_2q_3a_1^{n}a_2^{-1}q_4a_1^{-n}a_2^{-1}k_1q_5c^n\dots
k_Nq_5c^n=1$$ can be obtained by application of $O(n^2)$ relations
(7.1) -- (7.3) to the hub, and the relation $W_n(b_1,b_2)=1$ can
be deduced from (7.1) -- (7.10) in $O(n^2)$ steps.\end{lm}$\Box$

\medskip

The proof of the next claim is absolutely similar to the proof of Lemma \ref{3.1}.

\medskip

\begin{lm}\label{7.2}  A minimal hub-free diagram over the presentation
of the group $H_{k,1}$ has no $\rho$- $r$-, $(r,q)$-, $q$-, $(r,k)$-, $k$-,
$(\rho,k)$-, $(a,b)$-, $(\rho,a)$-, $(\rho,c)$-, $(r,a)$-, $(r,c)$-, $b$-,
$a$-, $(\rho,q)$- or $(q,b)$-annuli.\end{lm}  $\Box$

\medskip

The statement of Lemma \ref{4.1} is also true for the group $G=G_{k,1}$, because
the boundary label of a spoke has the form $r^\ell$, but in the HNN-extension
of a free group with the stable letter $r$ defined by relations (7.1)--(7.3),
the word $r^\ell$ can commute with a subword of the left-hand side of (7.4)
written between neighbor $k$-letters only if $\ell=0$. In view of minimality of $\Delta$,
this means that the spokes have length 0, and the two hubs form a mirror
pair of cells, that cancel.

\medskip

\begin{lm}\label{7.3} The statement of Lemma 4.2 is true for
diagrams over $H=H_{k,1}$. Hence Lemmas 4.3, 4.4 are also valid for $H$.
\end{lm}
\medskip

$\Box$ As in the proof of Lemma \ref{4.2}, we have to consider the boundary
label $V=V(\rho^{\pm 1},r^{\pm 1})$ of the spokes $T_1, T_2$ that
must commute with $q_5$ in view of relations (7.1)--(7.3), (7.4)--(7.11).
This is possible if and only if the exponent sum for $r$ is zero in
$V$, since the group defined by (7.1)--(7.3), (7.4)--(7.11) has
the retraction preserving $r, q_5$ and $c$ and mapping the other
generators to 1.

As in Lemma \ref{4.2} we have to prove that $V$ commutes with the cyclic
permutation $W$ of the left-hand side of (7.4) beginning with $k_j$.
The word $V$ commutes with any letter $k_j$ by relations (7.2), (7.8).
Therefore it suffices to prove that $V$ commutes with
the word $\Lambda_0\equiv q_1a_2q_2a_2q_3a_2^{-1}q_4a_2^{-1}$.

It is clear that the equality $r\Lambda_nr^{-1}=\Lambda_{n+1}$
follows from (7.1)--(7.3) for \\ $\Lambda_s\equiv
q_1a_1^sa_2q_2a_1^{-s}a_2q_3a_1^sa_2^{-1}q_4a_1^{-s}a_2^{-1}$ and
any integer $n$. The equalities $\rho
\Lambda_s\rho^{-1}=\Lambda_s$ follows from (7.5), (7.7),
(7.9)--(7.10) and Lemma \ref{7.1}. Since the exponent sum for the
occurrences of $r$ in $V$ is equal to 0, $V\Lambda_0V^{-1}=V_0$,
as desired. $\Box$

\medskip

Let us mention another similarity of the groups $H(v,m)$ and
$H_{k,1}$: the claims of Lemmas \ref{5.1}--\ref{5.4}
are true for $H=H_{k,1}$ as well,
and the proofs are quite analogous to those in Section 5.

To complete the proof of Theorem 2 we need

\medskip

\begin{lm} \label{7.4} Let $w=w(b_1^{\pm 1},b_2^{\pm 1})$ be a word of
length $n$ such that $w=1$ in $H$. Then this equality can be derived
from the trivial one by application of $O(n^4)$ of relations (7.1)--(7.11).\end{lm}

\medskip

$\Box$ Since $w=1$ in $BS_{k,1}$, the exponential sum over the
occurrences of $b_1$ in $w$ is equal to $0$. Therefore the word
$w$ is freely equal to a product $\prod_i v_i$ where $v_i\equiv
b_1^{s_i} b_2^{\pm 1}b_1^{-s_i}$, $|s_i| < n/2$ and the number of
factor do not exceed $n$. Passing to a freely conjugate word, we
may assume that $0\le s_i<n$ for every $i$.

Let $s=\min s_i$, and assume there are factors
$v_i\equiv b_1^sb_2b_1^{-s}$ and $v_j\equiv b_1^sb_2^{-1}b_1^{-s}$. Then by
Lemma \ref{7.1} we can transpose $v_i$ with a neighbor factor $v_{i\pm1}$ by applying the defining relations at most $O(n^2)$ times, and so after $O(n^3)$ applications
of the relations the factors $v_i$ and $v_j$ cancel.

Now assume that all the $v_i$'s with $s_i=s$ are equal, and $\ell$ is the number
of such $v_i$'s.
Since $b_1^tb_2b_1^{-t}=b_2^{k^t}$ for any $t\ge s$ in $G$, and the product
of the words $v_i$ is equal to 1 in $G$, the number $\ell$ is a multiple of $k$.

Again, by Lemma \ref{7.1} we can collect the $\ell=km$
factors (with minimal $s_i$'s) at the end of the word $w$ applying
the relations at most $\ell O(n^3)$ times. This suffix is (in the free group)
the product of $m$ factors $u_i\equiv b_1^sb_2^{\pm k}b_1^{-s}$, and
applying relator (7.11) $m$ times we can rewrite it as the product of
$m$ factors $b_1^{s+1}b_2^{\pm 1}b_1^{-s-1}$.

Thus, we need $O(n^3)$ relations to reduce the word $w$ to 1, or $\ell O(n^3)$
relations to decrease the number of factors $v_i$ by $\ell(k-1)$.
The lemma is proved. $\Box$

\medskip

Now an isoperimetric inequality for the group $H=H_{k,1}$ can be obtained by
almost the same reasoning as in Lemma \ref{6.1} for $H=H(v,m)$ with the only
essential difference that the function $f(n^2)^2$ should be replaced by $(n^2)^4$,
since by Lemma 7.4 the ``area" of a $B$-cell of perimeter $s$ (in the
relators (7.1)--(7.11)) is at most $O(s^4)$. The property that the constructed
in the proof of Theorem 2 embedding is undistorted, can be justified in the
same manner as for Theorem 1.

\bigskip
\noindent Alexander Yu. Olshanskii\\
Department of Mathematics\\
Vanderbilt University \\
olsh@math.vanderbilt.edu\\
and\\
Department of Higher Algebra\\
MEHMAT\\
 Moscow State University\\
olsh@nw.math.msu.su\\

\noindent Mark V. Sapir\\
Department of Mathematics\\
Vanderbilt University\\
http://www.math.vanderbilt.edu/$\sim$msapir\\

\end{document}